\documentclass[11pt]{amsart}
\usepackage{amssymb,psfig}

\DeclareMathOperator{\dom}{Dom}
\newcommand{\dbar}{\ensuremath{\overline\partial}}
\newcommand{\dbarstar}{\ensuremath{\overline\partial^*}}
\newcommand{\C}{\ensuremath{\mathbb{C}}}

\makeatletter
\newcommand{\sumprime}{\if@display\sideset{}{'}\sum%
            \else\sum'\fi}
\makeatother

\begin{document}

\numberwithin{equation}{section}

\newtheorem{theorem}{Theorem}[section]
\newtheorem{proposition}[theorem]{Proposition}
\newtheorem{conjecture}[theorem]{Conjecture}
\def\theconjecture{\unskip}
\newtheorem{corollary}[theorem]{Corollary}
\newtheorem{lemma}[theorem]{Lemma}
\newtheorem{observation}[theorem]{Observation}
\theoremstyle{definition}
\newtheorem{definition}{Definition}
\numberwithin{definition}{section}
\newtheorem{remark}{Remark}
\def\theremark{\unskip}
\newtheorem{question}{Question}
\def\thequestion{\unskip}
\newtheorem{example}{Example}
\def\theexample{\unskip}
\newtheorem{problem}{Problem}

\def\intprod{\mathbin{\lr54}}
\def\reals{{\mathbb R}}
\def\integers{{\mathbb Z}}
\def\N{{\mathbb N}}
\def\complex{{\mathbb C}\/}
\def\distance{\operatorname{distance}\,}
\def\spec{\operatorname{spec}\,}
\def\interior{\operatorname{int}\,}
\def\cl{\operatorname{cl}\,}
\def\essspec{\operatorname{esspec}\,}
\def\range{\operatorname{\mathcal R}\,}
\def\kernel{\operatorname{\mathcal N}\,}
\def\linearspan{\operatorname{span}\,}
\def\lip{\operatorname{Lip}\,}
\def\sgn{\operatorname{sgn}\,}
\def\Z{ {\mathbb Z} }
\def\e{\varepsilon}
\def\p{\partial}
\def\rp{{ ^{-1} }}
\def\Re{\operatorname{Re\,} }
\def\Im{\operatorname{Im\,} }
\def\dbarb{\bar\partial_b}
\def\eps{\varepsilon}

\def\Hs{{\mathcal H}}
\def\E{{\mathcal E}}
\def\scriptu{{\mathcal U}}
\def\scriptr{{\mathcal R}}
\def\scripta{{\mathcal A}}
\def\scripti{{\mathcal I}}
\def\scriptk{{\mathcal K}}
\def\scripth{{\mathcal H}}
\def\scriptm{{\mathcal M}}
\def\scripte{{\mathcal E}}
\def\scriptt{{\mathcal T}}
\def\scriptb{{\mathcal B}}
\def\scriptf{{\mathcal F}}
\def\scripto{{\mathfrak o}}
\def\scriptv{{\mathcal V}}
\def\frakg{{\mathfrak g}}
\def\frakG{{\mathfrak G}}

\def\ov{\overline}
\date {February, 2004.}

\author{Siqi Fu}
\thanks
{The author was supported in part by NSF grant DMS 0070697 and by
an AMS centennial fellowship.}
\address{Department of Mathematical Sciences,
Rutgers University-Camden, Camden, NJ 08102}
\email{sfu@camden.rutgers.edu}
\title[]  
{Hearing pseudoconvexity with the Kohn Laplacian} \maketitle



\section{Introduction}

Mark Kac's famous question ``Can one hear the shape of a drum?"
asks whether the spectrum of the Dirichlet Laplacian determines a
planar domain up to congruence~\cite{Kac66}.  This question was
answered negatively by Gordon, Webb, and
Wolpert(cf.~\cite{GordonWebbWolpert92}).  It has inspired a
tremendous amount of research on the interplay of the spectrum of
differential operators and the geometry of  ambient spaces. Here
we study the several complex variables analogue of Kac's question:
To what extent is the geometry of a bounded domain $\Omega$ in
$\C^n$ determined by the spectrum of the $\dbar$-Neumann and Kohn
Laplacians? Since the work of Kohn~\cite{Kohn63, Kohn64}, it has
been discovered that various notions of regularity of the
$\dbar$-Neumann and Kohn Laplacians, such as subellipticity,
hypoellipticity, and compactness, are intimately related to the
boundary geometry of the domain. (See, for example, the surveys
\cite{BoasStraube99, Christ99, DangeloKohn99, FuStraube01}.) It is
then natural to expect that one should be able to ``hear" more
about the geometry of a bounded domain in $\C^n$ with the
$\dbar$-Neumann and Kohn Laplacians than with the usual Dirichlet
Laplacians. In this paper, we prove the following:

\begin{theorem}\label{maintheorem}
Let $\Omega$ be a bounded domain in $\C^n$, $n>1$, with connected
Lipschitz boundary $b\Omega$.  Let $\square_{b, q}$ be the Kohn
Laplacian on $L^2_{(0, q)}(b\Omega)$. Let
$\essspec(\square_{b,q})$ be the essential spectrum of
$\square_{b, q}$.  If $\inf\essspec(\square_{b, q})>0$ for all
$1\le q\le n-1$, then $\Omega$ is pseudoconvex.
\end{theorem}

It was shown by Kohn~\cite{Kohn86} that on smooth pseudoconvex
boundaries $b\Omega$ in Stein manifolds, $\dbar_b$ has closed
range in $L^2_{(0, q)}(b\Omega)$ for all $1\le q\le n-1$.
Independently, Shaw~\cite{Shaw85} (for $1\le q\le n-2$) and
Boas-Shaw~\cite{BoasShaw86} (for $q=n-1$) established
$L^2$-existence theorems for the $\dbar_b$-operator on smooth
pseudoconvex boundaries in $\C^n$. Recently, Shaw~\cite{Shaw03}
extended these results to pseudoconvex Lipschitz boundaries.  In
light of these results and Theorem~\ref{maintheorem}, for
connected and sufficiently smooth boundaries in $\C^n$,
pseudoconvexity is characterized by positivity of the infimum of
the spectrum (or the essential spectrum) of the Kohn Laplacians on
all $(0, q)$-forms, $1\le q\le n-1$.

This paper is organized as follows.  In Section~\ref{pre}, we
recall necessary setups and definitions. Section~\ref{psc}
contains the proof of Theorem~\ref{maintheorem}. Further remarks
are given in Section~\ref{remarks}.

\bigskip

\noindent{\bf Acknowledgment:} Part of this work was done while
the author visited Princeton University on an AMS Centennial
Research Fellowship. The author is indebted to Professors J. J.
Kohn, M.-C. Shaw, and Y.-T. Siu for helpful conversations and kind
encouragement.

\section{Preliminaries}\label{pre}

We first review the well-known operator theoretic setup
(cf.~\cite{Hormander65, FollandKohn72}).  Let $T_k\colon H_k\to
H_{k+1}$, $k=1, 2$, be densely defined, closed operators between
Hilbert spaces. Assume that $\range(T_1)\subset \kernel(T_2)$,
where $\range$ and $\kernel$ denote the range and kernel of the
operators. Let $T^*_k$ be the Hilbert space adjoint of $T_k$. Then
$T^*_k$ is also densely defined and closed.  Let
\[
Q(u, v)=(T^*_1 u, T^*_1 v)+(T_2 u, T_2 v)
\]
with $\dom(Q)=\dom(T^*_1)\cap\dom(T_2)$.  It is easy to see that
$Q(u, v)$ is a non-negative, densely defined, closed sesquilinear
form on $H_2$.  It follows that $Q(u, v)$ uniquely determines a
non-negative, densely defined, self-adjoint operator $\square$ on
$H_2$ such that $\dom(\square^{1/2})=\dom(Q)$ and $Q(u,
v)=(\square u, v)$ for all $u\in\dom(\square)$ and $v\in\dom(Q)$.
(We refer the reader to \cite{Davies95, Kato76, ReedSimon} for
detail on sesquilinear forms and self-adjoint operators.) The
spectrum $\spec(\square)$ of $\square$ is a non-empty closed
subset of $[0, \infty)$ and the infimum of the spectrum is given
by
\[
\inf\spec (\square) =\inf\{ Q(u, u); \ \ u\in\dom(Q), \| u\|=1 \}.
\]
For any positive integer $j$, let
\[
\lambda_j=\sup_{v_1, \ldots, v_{j-1}\in \dom(Q)}\inf\left\{ Q(u,
u); \ u\in\dom(Q), u\perp v_i, 1\le i\le j-1, \|u\|=1\right\}.
\]
Then $\square$ has compact resolvent if and only if
$\lambda_j\to\infty$. In this case, $\lambda_j$ is the $j^{\rm
th}$ eigenvalue of $\square$, when the eigenvalues are arranged in
increasing order and repeated according to multiplicity.  If
$\square$ has non-compact resolvent (equivalently, the essential
spectrum $\essspec(\square)$ is non-empty), $\lambda_j$ is either
an eigenvalue of finite multiplicity or the bottom of
$\essspec(\square)$.  In either cases,
$\lim_{j\to\infty}\lambda_j=\inf\essspec(\square)$. In what
follows, we will set $\inf\essspec(\square)=\infty$ when
$\essspec(\square)$ is empty.

\begin{lemma}\label{spectral lemma} With the above notations
and assumptions, $\inf\spec(\square)>0$ if and only if
$\range(T_2)$ is closed and $\range(T_1)=\kernel(T_2)$.
Furthermore, $\inf\essspec(\square)>0$ if and only if there exists
a finite dimensional subspace $L\subset\dom(Q)$ such that
$\range(T_2\big|_{L^\perp})$ is closed and $\range(T_1)\cap
L^\perp=\kernel(T_2)\cap L^\perp $.
\end{lemma}

The first part of the lemma is well-known (compare
\cite{Hormander65}, Theorem 1.1.2; \cite{Catlin83}, Proposition 3;
and \cite{Shaw92}, Proposition 2.3). We provide a proof here for
completeness. To prove the forward direction, we note that
$\inf\spec(\square)>0$ implies that $\square$ has a bounded
inverse $N$ defined on all $H_2$. Hence each $u\in H_2$ has an
orthogonal decomposition $u=T_1 T_1^*Nu +T_2^*T_2 Nu$. It follows
that $\range(T_1)=\kernel(T_2)$ and
$\range(T_2^*)=\kernel(T_1^*)$. Since now $T_2^*$ has closed
range, so is $T_2$. We thus conclude the prove of forward
direction.  To prove the opposite, for any $u\in\dom(Q)$, we write
$u=u_1+u_2$ where $u_1\in \dom(Q)\cap\kernel(T_2)$ and
$u_2=\dom(Q)\cap\kernel(T_2)^\perp$. Since
$\kernel(T_2)=\range(T_1)=\kernel(T_1^*)^\perp$ and
$\kernel(T_2)^\perp=\range(T_1)^\perp=\kernel(T_1^*)$, there
exists a positive constant $C$ such that
$\|u\|^2=\|u_1\|^2+\|u_2\|^2\le C(\|T_1^* u_1\|^2+\|T_2 u_2\|^2)=C
Q(u, u)$. This concludes the proof of the backward direction.

For a proof of the second part of the lemma, we observe that by
the above-mentioned spectral theoretic results,
$\inf\essspec(\square)>0$ if and only if there exists a positive
constant $C$ and a finite dimensional subspace $L$ of $\dom(Q)$
such that
\[
Q(u, u)\ge C\|u\|, \quad u\in\dom(Q)\cap L^\perp .
\]
To prove the forward direction, let $H_2'=H_2\ominus L$ and let
$T_2'=T_2\big|_{H_2'}$ and ${T_1^*}'=T_1^*\big|_{H_2'}$.  Then
$T_2'\colon H_2'\to H_3$ and ${T_1^*}'\colon H_2'\to H_1$ are
densely defined, closed operators. Let $T_1'\colon H_1\to H_2'$ be
the adjoint of ${T_1^*}'$.  It is easy to see that
$\range(T_1')\subset\kernel(T_2')$ and $\dom(T_1')=\dom(T_1)$.
Applying the first part of the lemma to the operators $T'_1\colon
H_1\to H'_2$ and $T'_2\colon H'_2\to H_3$ and the sesquilinear
form
\[
Q'(u, v)=({T_1'}^* u, {T_1'}^*v)+(T_2' u, T_2'v)
\]
on $H_2'$ with $\dom(Q')=\dom(Q)\cap L^\perp$, we obtain that
$T'_1$ and $T'_2$ have closed range and
$\range(T'_1)=\kernel(T'_2)$. We then conclude the proof of the
forward direction by noting that $\range(T'_1)=\range(T_1)\cap
L^\perp$ and $\kernel(T'_2)=\kernel(T_2)\cap L^\perp$.  The
converse is treated similarly as above and is left to the reader.

\smallskip
\noindent{\bf Remark}.  Let
$\widetilde{H_2}=\kernel(T_1^*)^\perp$. Let $\widetilde
T_1^*=T_1^*\big|_{\widetilde{H}_2}$ and let $\widetilde Q(u,
v)=(\widetilde T_1^*u,\ \widetilde T_1^* v)$ be the sesquilinear
form on $\widetilde{H}_2$ with
$\dom(\widetilde{Q})=\dom(T_1^*)\cap\widetilde{H}_2$.  Let
$\widetilde\square$ be the self-adjoint operator determined by
$\widetilde{Q}(u, v)$. In this case,
$\inf\spec(\widetilde{\square})>0$ if and only if
$\range(T_1)=\kernel(T_1^*)^\perp$, and
$\inf\essspec(\widetilde{\square})>0$ if and only if there exists
a finite dimensional subspace $L$ of $\widetilde{H}_2$ such that
$\range(T_1)\cap L^\perp=L^\perp$.

\smallskip

We now review the $\dbar_b$-complex as introduced by
Kohn~\cite{Kohn65, KohnRossi65}, and adapted to Lipschitz
boundaries by Shaw~\cite{Shaw03}.  Let $\Omega$ be a bounded
Lipschitz domain in $\C^n$. (Recall that $b\Omega$ is Lipschitz if
it is given locally by a Lipschitz graph.) Let $\rho\in
\lip(\C^n)$ be a defining function of $b\Omega$ such that $\rho<0$
on $\Omega$ and $C_1\le |d\rho|\le C_2$ a.e. on $b\Omega$ for some
positive constants $C_1$ and $C_2$ (cf.~\cite{Shaw03}). Let $I^{0,
q}$, $0\le q\le n$, be the ideal in $\Lambda^{0, q} T^*(\C^n)$
generated by $\rho$ and $\dbar\rho$. Let $\Lambda^{0, q}
T^*(b\Omega)$ be the orthogonal complement with respect to the
standard Euclidean metric  of $I^{0, q}|_{b\Omega}$ in
$\Lambda^{0, q} T^*(\C^n)|_{b\Omega}$. Let $\tau\colon\Lambda^{0,
q} T^*(\C^n)|_{b\Omega}\to \Lambda^{0, q} T^*(b\Omega)$ be the
orthogonal projection.

Let $L^2_{(0, q)}(b\Omega)$ be the space of $(0, q)$-forms with
$L^2$-coefficients, equipped with the induced Euclidean metric on
$b\Omega$; that is, the projections under $\tau$ of $(0, q)$-forms
on $\C^n$ whose coefficients are in $L^2(b\Omega)$ when restricted
to $b\Omega$. The operator $\dbar_{b, q}\colon L^2_{(0,
q)}(b\Omega)\to L^2_{(0, q+1)}(b\Omega)$, $0\le q\le n-1$, defined
in the sense of distribution as the restriction of $\dbar_q$ to
the boundary $b\Omega$, is densely defined and closed
(see~\cite{Shaw03}). Let $\dbarstar_{b, q}$ be the Hilbert space
adjoint of $\dbar_{b, q}$. Let
\[
Q_{b, q}(u, v)=(\dbar_{b, q} u, \ \dbar_{b, q} v) +(\dbarstar_{b,
q-1} u, \ \dbarstar_{b, q-1} v)
\]
with $\dom(Q_{b, q})=\dom(\dbar_{b, q})\cap \dom(\dbarstar_{b,
q-1})$ when $1\le q\le n-2$, and let
\[
Q_{b, n-1}(u, v)=(\dbarstar_{b, n-2} u,\ \dbarstar_{b, n-2} v)
\]
with $\dom(Q_{b, n-1})=\dom(\dbarstar_{b, n-2})\cap
\kernel(\dbarstar_{b, n-2})^\perp$.  Then $Q_{b, q}$, $1\le q\le
n-1$, are non-negative, closed, and densely defined sesquilinear
forms on $L^2_{(0, q)}(b\Omega)$.  Therefore it uniquely
determines a non-negative, closed, densely defined, and
self-adjoint operator $\square_{b, q}$ on $L^2_{(0, q)}(b\Omega)$
such that $\dom(\square_{b, q}^{1/2})=\dom (Q_{b, q})$ and $Q_{b,
q}(u, v)= (\square_{b, q}u, v)$ for all $u\in\dom(\square_{b, q})$
and $v\in\dom(Q_q)$. The Kohn Laplacian is formally given by
$\square_{b, q}=\dbar_{b, q-1}\dbarstar_{b, q-1}+ \dbarstar_{b,
q}\dbar_{b, q}$ for $1\le q\le n-2$ and $\square_{b,
n-1}=\dbar_{b, n-2}\dbarstar_{b,
n-2}\big|_{\kernel(\dbarstar_{n-2})^\perp}$. (Notice that on top
degree $(0, n-1)$-forms, the Kohn Laplacian here is the
restriction to the orthogonal complement of $\kernel(\dbarstar_{b,
n-2})$ of the usual Kohn Laplacian.  We make this restriction
because the kernel of $\dbarstar_{b, n-2}$ is infinite
dimensional.)  We refer the reader to the monographs [FK72] and
[CS01] for detail on the subject.

\section{Proof of the Main Theorem}\label{psc}

Let $\rho\in \lip(\C^n)$ be a global defining function of $\Omega$
such that $\rho<0$ on $\Omega$ and  $C_1\le |d\rho|\le C_2$ a.e.
on $b\Omega$. Arguing via {\it reductio ad absurdum}, we assume
that $\Omega$ is not pseudoconvex. Then there exists a domain
$\widetilde\Omega\supsetneqq\Omega$ such that every holomorphic
function on $\Omega$ extends holomorphically to $\widetilde\Omega$
(cf.~\cite{Hormander91}). Since $b\Omega$ is Lipschitz,
$\widetilde\Omega\setminus\cl(\Omega)$ is non-empty. After a
translation and a unitary transformation, we may assume that the
origin is in $\widetilde\Omega\setminus\cl(\Omega)$ and the
$z_n$-axis has a non-empty intersection with $\Omega$.
Furthermore, we may assume that the positive $y_n$-direction is
the outward normal direction of the intersection of the $y_n$-axis
with $b\Omega$ and $b\Omega\cap\widetilde{\Omega}$ is
parameterized near the intersection by $y_n=h(z_1, \ldots,
z_{n-1}, x_n)$ for some Lipschitz function $h$.

For any integers $\alpha \ge 0$, $m\ge 1$, and $q\ge 1$,  and for
any $\{k_1, \ldots, k_{q-1}\}\subset \{1, 2, \ldots, n-1\}$, let
\[
u_{\alpha, m}(k_1, \ldots, k_q)=\frac{(\alpha+q-1)! \bar
z_n^{m\alpha}(\bar z_{k_1}\cdots\bar z_{k_q})^{m-1}}{ r_m^{\alpha+
q}}\sum_{j=1}^q (-1)^j \bar{z}_{k_j} d\bar{z}_{k_1}\wedge
\ldots\wedge\widehat{d\bar{z}_{k_j}}\wedge\ldots\wedge
d\bar{z}_{k_q}
\]
where $k_q=n$, $r_m=|z_1|^{2m}+\ldots +|z_n|^{2m}$, and
$\widehat{d\bar{z}_{k_j}}$ indicates as usual the omission of
$d\bar z_{k_j}$ from the wedge product. It is evident that
$u_{\alpha, m}(k_1, \ldots, k_q)$ is a smooth $(0, q-1)$-form on
$\C^n\setminus\{ 0\}$ that is skew-symmetric with respect to the
indices $(k_1, \ldots, k_{q-1})$.  In particular, $u_{\alpha,
m}(k_1, \ldots, k_q)=0$ when two $k_j$'s are identical. Write
$K=(k_1,\ldots, k_q)$, $d\bar z_K=d\bar z_{k_1}\wedge \ldots\wedge
d\bar z_{k_q}$, $\bar z_K^{m-1}= (\bar z_{k_1}\cdots\bar
z_{k_q})^{m-1}$, and $\widetilde{d\bar
z_{k_j}}=d\bar{z}_{k_1}\wedge
\ldots\wedge\widehat{d\bar{z}_{k_j}}\wedge\ldots\wedge
d\bar{z}_{k_q}$.  Then

\begin{align*}
\dbar u_{\alpha, m}(k_1, \ldots, k_q)&=-\frac{(\alpha+q)!m \bar
z_n^{m\alpha}\bar z_K^{m-1}} {r_m^{\alpha+q+1}}\big(r_m d\bar z_K
+ \big(\sum_{\ell =1}^n \bar z_\ell^{m-1} z_\ell^m d\bar
z_\ell\big)\wedge \big(\sum_{j=1}^q (-1)^j \bar
z_{k_j}\widetilde{d\bar z_{k_j}}\big)\big)\\
&=-\frac{(\alpha+q)!m\bar z_n^{m\alpha}\bar
z_K^{m-1}}{r_m^{\alpha+q+1}} \sum_{\ell\in\{1, \ldots,
n\}\setminus\{k_1, \ldots, k_q\}} z^m_\ell \bar
z_\ell^{m-1}\big(\bar z_\ell d\bar z_K+\sum_{j=1}^q(-1)^j \bar
z_{k_j}
\widetilde{d\bar z_{k_j}}\big)\\
&=m\sum_{\ell=1}^{n-1} z^m_\ell u_{\alpha, m}(\ell, k_1, \ldots,
k_q).
\end{align*}

In particular, $u_{\alpha, m}(1, \ldots, n)$ is $\dbar$-closed.
Let $N=(1/|\partial\rho|)\sum_{j=1}^n \rho_{z_j}
\partial/\partial\bar z_j$ and let
$$
u^b_{\alpha, m}(k_1, \ldots, k_q)=\tau(u_{\alpha, m}(1, 2, \ldots,
n)) =N\lrcorner(\frac{\dbar\rho}{|\dbar\rho|}\wedge u_{\alpha,
m}(k_1, \ldots, k_q))\in L^2_{(0, q-1)}(b\Omega),
$$
where $\lrcorner$ denotes the contraction operator.  Then for
$1\le q\le n-1$,
\[
\dbar_{b, q-1} u^b_{\alpha, m}(k_1, \ldots,
k_q)=m\sum_{\ell=1}^{n-1} z^m_\ell u^b_{\alpha, m}(\ell, k_1,
\ldots, k_q).
\]

We now show that  $u^b_{\alpha, m}(1, 2, \ldots, n)\perp
\kernel(\dbarstar_{b, n-1})$.   Let $\star\colon L^2_{(p,
q)}(\Omega) \to L^2_{(n-p, n-q)}(\Omega)$ be the Hodge star
operator, defined by $\langle \phi, \psi\rangle
dV=\phi\wedge\star\psi$ where $dV$ is the Euclidean volume form.
Let $v\in\kernel(\dbarstar_{b, n-1})$. Let
$\theta=\star(dz_1\wedge\ldots\wedge dz_n\wedge
\dbar\rho/|\dbar\rho|)$. Then $v=\bar f \theta$ for some $f\in
L^2(b\Omega)$ with $\dbar_b f=0$. It follows from a version of
Hartogs-Bochner extension theorem that there exists a holomorphic
function $F$ on $\Omega$ such that the non-tangential limit of $F$
agrees with $f$ a.e. on $b\Omega$, and
$$
\lim_{\epsilon\to 0^+} \int_{b\Omega} |F(z-\epsilon \nu(z)) -f(z)|^2
d\sigma=0
$$
where $\nu(z)=\nabla \rho/|\nabla \rho|$. (See, for example,
Theorem 7.1 in \cite{Kytmanov95}.  Although the theorem is stated
only for $C^1$-smooth boundaries, the proof works for Lipschitz
boundaries with only minor modifications.)  Let $\nu_\delta (z)$
be the convolution of $\nu(z)$ with appropriate Friederichs'
mollifiers. Then there exists a subsequence $\delta_j\to 0$ such
that $\nu_{\delta_j}(z)\to \nu(z)$ a.e. on $b\Omega$.  Therefore,
\begin{align*}
(u^b_{\alpha, m}(1, \ldots, n), v) &=\int_{b\Omega}
f(z)u^b_{\alpha, m}(1, \ldots, n)(z)
\wedge dz_1\ldots\wedge dz_n \\
&=\lim_{\epsilon\to 0} \int_{b\Omega} F(z-\epsilon \nu(z))
u^b_{\alpha, m}(1,
\ldots, n)(z)\wedge dz_1\ldots \wedge dz_n \\
&=\lim_{\epsilon\to 0}\lim_{\delta_j\to 0} \int_{b\Omega}
F(z-\epsilon \nu_{\delta_j}(z))
u^b_{\alpha, m}(1, \ldots, n)(z)\wedge dz_1\ldots \wedge dz_n\\
&=\lim_{\epsilon\to 0}\lim_{\delta_j\to 0}\int_\Omega \dbar\big(
F(z-\epsilon\nu_{\delta_j}(z)) u_{\alpha, m}(1, \ldots,
n)(z)\wedge dz_1\ldots\wedge dz_n\big) =0 .
\end{align*}
Hence $u^b_{\alpha, m}(1, \ldots, n)\perp \kernel(\dbarstar_{b,
n-1})$ as claimed.

By Lemma~\ref{spectral lemma} and the subsequence remark, we can
choose a sufficiently large positive integer $M$ such that there
exist subspaces $S_q$ of $\dom(Q_{b, q})$ for $1\le q\le n-2$ and
$S_{n-1}$ of $\kernel(\dbarstar_{b, n-2})^\perp$, all of which
have dimensions $<M$ and satisfy $\range(\dbar_{b, q-1})\cap
S_q^\perp=\kernel(\dbar_{b, q})\cap S_q^\perp$, $1\le q\le n-2$,
and $\range(\dbar_{b, n-2})\cap S_{n-1}^\perp=S_{n-1}^\perp$. Fix
$m\ge 1$ (to be specified later) and let $\scriptf_0$ be the
linear span of $\{u^b_{\alpha, m}(1, \ldots, n); \ \alpha=1,
\ldots, M^{n-1}\}$. For any $u\in\scriptf_0$ and for any $\{k_1,
\ldots, k_{q-1}\}\subset\{1, \ldots, n-1\}$, we set
\[
u(k_1, \ldots, k_{q-1}, n)=\sum_{j=1}^k c_j u^b_{\alpha_j, m}(k_1,
\ldots, k_{q-1}, n)
\]
if $u=\sum_{j=1}^k c_j u^b_{\alpha_j, m}(1, \ldots, n)$.  We
decompose $\scriptf_0$ into a direct sum of $M^{n-2}$ subspaces,
each of which is $M$-dimensional. Since $\dim(S_{n-1})<M$ and
$u_{\alpha, m}(1, \ldots, n)\in \kernel(\dbarstar_{b,
n-2})^\perp$, there exists a non-zero form $u$ in each of the
subspaces such that $\dbar_b v_u(\emptyset)=u$ for some
$v_u(\emptyset)\in L^2_{(0, n-2)}(b\Omega)$.  Let $\scriptf_1$ be
the $M^{n-2}$-dimensional linear span of all such $u$'s.  We
extend $u\mapsto v_u(\emptyset)$ linearly to all $u\in\scriptf_1$.

For $0\le q\le n-1$, we use induction on $q$ to construct an
$M^{n-q-2}$-dimensional subspace $\scriptf_{q+1}$ of $\scriptf_q$
with the properties that for any $u\in\scriptf_{q+1}$, there
exists $v_u(k_1, \ldots, k_q)\in L^2_{(0, n-q-2)}(b\Omega)$ for
all $\{k_1, \ldots, k_q\}\subset\{1, \ldots, n-1\}$ such that

\begin{enumerate}

\item $v_u(k_1, \ldots, k_q)$ depends linearly on $u$.

\item $v_u(k_1, \ldots, k_q)$ is skew-symmetric with respect to
indices $K=(k_1, \ldots, k_q)$.

\item $\dbar_b v_u(K)=m\sum_{j=1}^q (-1)^j z^m_{k_j} v_u(K;
\hat{k_j}) + (-1)^{q+|K|} u(1, \ldots, n; \hat K) $ where
$|K|=k_1+ \ldots +k_q$.  The hat $\hat{ }$\; indicates deletion of
indices beneath it from the indices preceding the semicolon in the
same enclosing parenthesis.
\end{enumerate}

We now show how to construct $\scriptf_{q+1}$ and $v_u(k_1,
\ldots, k_q)$ for $u\in\scriptf_{q+1}$ and $\{k_1,\ldots,
k_q\}\subset\{1, \ldots, n-1\}$ once $\scriptf_q$ has been
constructed.  For any $u\in\scriptf_q$ and any $\{k_1, \ldots,
k_q\} \subset\{1, \ldots, n-1\}$, write $K=(k_1, \ldots, k_q)$,
and let
\[
w_u(K)=m\sum_{j=1}^q (-1)^j z^m_{k_j} v_u(K; \hat{k_j})
+(-1)^{q+|K|} u(1, \ldots, n; \hat K) .
\]
Then
\begin{align*}
\dbar_b w_u(K)&=m\sum_{j=1}^q (-1)^j z^m_{k_j}\dbar_b v_u (K;
\hat{k_j})
+(-1)^{q+|K|} \dbar_b u(1, \ldots, n; \hat{K})\\
&=m\sum_{j=1}^q (-1)^j z^m_{k_j} \big( m\sum_{1\le i <j} (-1)^i
z^m_{k_i} v_u(K; \hat k_j, \hat k_i) +m\sum_{j<i\le q} (-1)^{i-1}
z^m_{k_i}  v_u(K; \hat k_j,
\hat k_i) \\
&\qquad -(-1)^{q+|K|-k_j} u(1, \ldots, n; \widehat{(K; \hat k_j)})
\big) +(-1)^{q +|K|} \dbar_b u(1, \ldots, n; \hat K) \\
&=(-1)^{q+|K|}\big( -m\sum_{j=1}^q (-1)^{j-k_j} z^m_{k_j} u(1,
\ldots, n; \widehat{(K; \hat k_j)}) +\dbar_b
u(1, \ldots, n; \hat K)\big) \\
&=(-1)^{q+|K|}\big( -m\sum_{j=1}^q z^m_{k_j} u(k_j, (1, \ldots, n;
\hat K)) + \dbar_b u(1, \ldots, n; \hat K)\big)=0.
\end{align*}
We again decompose $\scriptf_q$ into a direct sum of $M^{n-q-2}$
linear subspaces, each of which is $M$-dimensional.  Since
$\dim(S_{n-q-2})<M$ and $\dbar_b w_u(K)=0$, there exists a
non-zero form $u$ in each of these subspaces such that $\dbar_b
v_u(K)=w_u(K)$ for some $v_u(K)\in L^2_{(0, n-q-2)}(b\Omega)$.
Since $w_u(K)$ is skew-symmetric with respect to indices $K$, we
may choose $v_u(K)$ to be skew-symmetric with respect to $K$ as
well.  The subspace $\scriptf_{q+1}$ of $\scriptf_q$ is then the
linear span of all such $u$'s.

Note that $\dim(\scriptf_{n-1})=1$.  Let $u$ be any non-zero form
in $\scriptf_{n-1}$ and let
\[
g=w_u(1, \ldots, n-1)= m\sum_{j=1}^{n-1} z^m_j v_u(1, \ldots, \hat
j, \ldots, n-1) - (-1)^{n+\frac{n(n-1)}2} u(n) .
\]
Then $g\in L^2(b\Omega)$ and $\dbar_b g =0$.  Therefore, $g$ has a
holomorphic extension $G$ to $\Omega$ such that the non-tangential
limit of $G$ agrees with $g$ a.e. on $b\Omega$ (cf. Theorem 7.1 in
\cite{Kytmanov95}). By the {\it reductio ad absurdum} assumption,
$G$ extends holomorphically to $\widetilde{\Omega}$.  Write
$z'=(z_1, \ldots, z_{n-1})$.  For sufficiently small $\e>0$ and
$\delta>0$,
\begin{align*}
&\int_{|x_n|<\e,|z'|<\e}
\big|\big(G+(-1)^{n+\frac{n(n-1)}2}u(n)\big)(\delta z',
x_n+ih(\delta z',
x_n))\big|dV(z')dx_n \\
&\qquad\le m\delta^m \sum_{j=1}^{n-1}\int_{|x_n|<\e, |z'|<\e}
|z_j|^m|v_u(1, \ldots, \hat j, \ldots, n-1)(\delta z',
x_n+ih(\delta z', x_n))|
dV(z')dx_n \\
&\qquad\le m\delta^{m-2(n-1)}\e^m\sum_{j=1}^{n-1} \|v_u(1, \ldots,
\hat j, \ldots, n-1)\|_{L^1(b\Omega)}.
\end{align*}
Choosing $m>2(n-1)$ and letting $\delta\to 0$, we obtain
\[
G(0, x_n+ih(0, x_n))=-(-1)^{n+\frac{n(n-1)}2} u(n)(0, x_n+ih(0,
x_n)).
\]
However, $u(n)(0, z_n)$ is a non-trivial linear combination of
functions of form $1/z^k$ with $k$ a positive integer. This leads
to a contradiction with the analyticity of $G$ near the origin. We
therefore conclude the proof of Theorem~\ref{maintheorem}.

\section{Further Remarks}\label{remarks}

(1) The analogue of Theorem~\ref{maintheorem} for the
$\dbar$-Neumann Laplacian $\square_q$ also holds under the
assumption that $\interior(\cl(\Omega))=\Omega$.  This is a
consequence of the sheaf cohomology theory (see \cite{Serre53,
Laufer66, Ohsawa88}), in light of Lemma~\ref{spectral lemma}. (We
thank Professor Y.-T. Siu for drawing our attention to
\cite{Laufer66}, by which the construction here is inspired.) The
above proof of Theorem~\ref{maintheorem} can be easily modified to
give a proof of this $\dbar$-Neumann Laplacian analogue, bypassing
sheaf cohomology arguments. In this case, one can actually choose
$m$ to be any positive integer, independent of the dimension $n$.
The non-elliptic nature of $\dbar_b$-complex seems to require that
the $m$ in the above proof be dependent on $n$. It follows from
H\"{o}rmander's $L^2$-existence theorem for the $\dbar$-operator
that $\inf\spec(\square_q)>0$ for all $1\le q\le n-1$ for any
bounded pseudoconvex domain in $\C^n$ (see \cite{Hormander65,
Hormander91}). Therefore, for a bounded domain $\Omega$ in $\C^n$
such that $\interior(\cl(\Omega))=\Omega$, the following
statements are equivalent: (a) $\Omega$ is pseudoconvex; (b)
$\inf\spec(\square_q)>0$ for all $1\le q\le n-1$; (c)
$\inf\essspec(\square_q)>0$ for all $1\le q\le n-1$.

(2) Let $\Omega$ be a bounded Lipschitz domain in $\C^n$ and let
$p\ge 1$. Consider $\dbar_{b, q}\colon L^p_{(0, q)}(\Omega)\to
L^p_{(0, q+1)}(b\Omega)$, $0\le q\le n-2$, where $L^p_{(0,
q)}(b\Omega)$ are boundary $(0, q)$-forms with $L^p$-coefficients.
Let $\scriptk_{n-1}$ be the space of all $f\in\dom(\dbar_{b,
n-1})$ such that
\[
\int_{b\Omega} f\wedge \alpha =0
\]
for all $\alpha\in C^\infty_{(n,
0)}(\ov{\Omega})\cap\kernel(\dbar)$.   Let
$H^p_q(b\Omega)=\kernel(\dbar_{b, q})/\range(\dbar_{b, q-1})$,
$1\le q\le n-2$, and
$H^p_{n-1}(b\Omega)=\scriptk_{n-1}/\range(\dbar_{b, n-1})$.  Then
the proof of Theorem~\ref{maintheorem} implies that $\Omega$ is
pseudoconvex if $\dim(H^p_q(b\Omega))<\infty$ for all $1\le q\le
n-1$.

(3) The generalization to $(p, q)$-forms is trivial.  We deal with
$(0, q)$-forms only for economy of notations.

\bibliography{survey}
%

\providecommand{\bysame}{\leavevmode\hbox to3em{\hrulefill}\thinspace}

\end{document}